\documentclass[reqno,11pt]{amsproc}
\usepackage{amssymb}
\usepackage{amscd}
\usepackage{amsmath}

\usepackage{xcolor}

\begin{document}

\title[Dynamic inverse problem for the one-dimensional system with memory.]
{Dynamic inverse problem for the one-dimensional system with
memory.}

\author[A.\,E. Choque-Rivero, A.\,S.~Mikhaylov, V.\,S.~Mikhaylov]
{$^1$A.\,E. Choque-Rivero, $^{2,3}$A.\,S.~Mikhaylov, $^2$V.\,S.~Mikhaylov}

\address{$^1$ Universidad Michoacana de San Nicolas de Hidalgo, Calle de Santiago Tapia 403, Centro, 58000 Morelia, Mexico,
    $^2$ St. Petersburg Department of V.A. Steklov Institute of
    Mathematics of the Russian Academy of Sciences, 7, Fontanka,
    191023 St. Petersburg, Russia. $^3$  St.Petersburg State University, 7/9 Universitetskaya
    nab., 199034 St. Petersburg,  Russia}

\email{mikhaylov@pdmi.ras.ru, vsmikhaylov@pdmi.ras.ru, abdon.ifm@gmail.com}

\begin{abstract}
We study the inverse dynamic  problem of recoverying the potential in the one-dimensional dynamical system with memory. The Gelfand--Levitan equations are derived for the kernel of the integral operator which is inverse to the control operator of the system. The potential is reconstructed  from the solution of these equations.

\end{abstract}


\maketitle

\newtheorem{corollary}{Corollary}
\newtheorem{definition}{Definition}
\newtheorem{lemma}{Lemma}
\newtheorem{proposition}{Proposition}
\newtheorem{remark}{Remark}
\newtheorem{theorem}{Theorem}

\section{Introduction}

For the function $K\in C^2_{loc}(\mathbb{R}_+)$ called the
relaxation kernel and the potential $q\in C^2_{loc}(\mathbb{R}_+)$
we consider the initial boundary value problem for the following
dynamical system:
\begin{equation}
\label{dyn_syst1}
\begin{cases}
u_{tt}(x,t)-u_{xx}(x,t)+q(x)u(x,t)+\int_0^tK(t-s)u(x,s)\,ds=0,\,\,
x> 0,\, t>0,\\
u(x,0)=0,\quad u_t(x,0)=0,\quad x\geqslant 0,\\
u(0,t)=f(t),\quad t\geqslant 0.
\end{cases}
\end{equation}
Here the function $f\in L_2(\mathbb{R}_+)$ is interpreted as a
\emph{boundary control}. The solution to (\ref{dyn_syst1}) is
denoted by $u^f$. The input-output correspondence in the system
(\ref{dyn_syst1}) is realized by the \emph{response operator} $R$
acting in $L_{2,\,loc}(\mathbb{R}_+)$ and defined by the formula
\begin{equation}
Rf=u^f_x(0,t).
\end{equation}
The finiteness of the speed of a signal propagating in the system
(\ref{dyn_syst1}) implies the following set up of the inverse
problem: assuming that the relaxation kernel is known and some
$T>0$ is fixed, one need to recover the potential $q(x)$, $x\in
(0,T)$ from the knowledge of the response operator
$R^{2T}:=R|_{L_2(0,2T)}$.

In \cite{Pan1,Pan2} the author considered the inverse dynamic
problem for the following dynamical system:
\begin{equation}
\label{dyn_syst}
\begin{cases}
v_t(x,t)-\int_0^tN(t-s)\left(v_{xx}(x,s)+\widetilde
q(x)v(x,s)\right)\,ds=0,\quad
0<x<L,\,t>0,\\
v(x,0)=0,\quad x\geqslant 0\\
v(0,t)=g(t),\quad v(L,t)=0\quad t\geqslant 0.
\end{cases}
\end{equation}
where the relaxation kernel $N\in C^3(0,L)$ such that $N(0)=1$,
the potential $\widetilde q\in C(0,L)$ and $g\in
L_{2,\,loc}(\mathbb{R}_+)$ is a boundary control. The solution of
(\ref{dyn_syst}) is denoted by $v^g(x,t)$. The input-output
correspondence for the system (\ref{dyn_syst}) is given by the
response operator defined by
\begin{equation}
\widetilde Rf=v^g_x(0,t),\quad t\geqslant 0.
\end{equation}
For this system the author studied the inverse problem of
recovering the potential $\widetilde q(x),$ $x\in (0,L)$ from the
knowledge of $\widetilde R|_{L_2(0,2L)}$. Note that the systems
(\ref{dyn_syst1}) and (\ref{dyn_syst}) are connected by special
change of unknown functions, known as "MacCamy trick", the details
are provided in \cite{Pan1}.

In \cite{Pan1,Pan2} the author partially used spectral methods in
solving the inverse problem, this, in particular requires the
boundary condition at $x=L$. We note that these methods are not
always adequate to the dynamic problems, and in our approach we
use purely dynamic tools.

In the second section we derive the Duhamel type representation
for the solution to (\ref{dyn_syst1}) and introduce operators of
the Boundary Control method \cite{B07,B17}. In the last section we
derive the Gelfand-Levitan equations and recover the unknown
potential.

\section{Forward problem. Operators of the Boundary control method.}

First we derive the Duhamel type representation for the function $u^f$:
\begin{lemma}
The solution to (\ref{dyn_syst1}) admits the following
representation
\begin{equation}
\label{u_repr} u^f(x,t)=f(t-x)+\int_x^tw(x,s)f(t-s)\,ds,
\end{equation}
where $w(x,s)$ is a solution to the following Goursat problem:
\begin{equation}
\label{Goursat}
\begin{cases}
w_{tt}-w_{xx}+qw+\int_0^tK(t-s)w(x,s)\,ds+K(t-x)=0,\ x > 0, t > 0,\\
\frac{d}{dx}w(x,x)=-\frac{1}{2}q(x),\quad x\geqslant 0,\\
w(x,0)=0,\quad x\geqslant 0.
\end{cases}
\end{equation}

\end{lemma}
\begin{proof}
The representation (\ref{u_repr}) implies that
\begin{align}
u^f_{tt}(x,t)=&f''(t-x)+\int_x^tw_{ss}(x,s)f(t-s)\,ds \nonumber\\
&+w_s(x,x)f(t-x)+w(x,x)f'(t-x),\label{e1}\\
u^f_{xx}(x,t)=&f''(t-x)-\frac{d}{dx}w(x,x)f(t-x)+w(x,x)f'(t-x)\nonumber\\
&-w_x(x,x)f(t-x)+\int_x^tw_{xx}(x,s)f(t-s)\,ds\label{e2}
\end{align}
Plugging these expressions into (\ref{dyn_syst1}) yields the relation
\begin{align}
0=&\int_x^t\left(w_{ss}(x,s)-w_{xx}(x,s)+q(x)w(x,s)\right)f(t-s)\,ds\nonumber\\
&+\left(2\frac{d}{dx}w(x,x)+q(x)\right)f(t-x)\nonumber\\
&+\int_0^tK(t-s)f(s-x)\,ds+\int_0^tK(t-s)\int_x^sw(x,\tau)f(s-\tau)\,d\tau\,ds\label{e3}
\end{align}
The first term in (\ref{e3}) after the change of variables
$s-x=t-\tau$ turns into

\begin{equation}
\int_x^{x+t}K(s-x)f(t-s)\,ds.\label{e4}
\end{equation}

The second term in (\ref{e3}) after the change of variables
$\tau=s-\alpha$ can be rewritten as
\begin{equation*}
\int_0^tK(t-s)\int_0^{s-x}w(x,s-\alpha)f(\alpha)\,d\alpha\,ds
\end{equation*}
after the changing the order of the integration we get
\begin{multline}
-\int_{-x}^0\int_0^{x+\alpha}K(t-s)w(x,s-\alpha)f(\alpha)dsd\alpha+\\
 \int_0^{t-x}\int_{\alpha+x}^tK(t-s)w(x,s-\alpha)f(\alpha)\,ds\,d\alpha=[s-\alpha=l]\\
-\int_{-x}^0\int_0^{x+\alpha}K(t-s)w(x,s-\alpha)f(\alpha)dsd\alpha+\\
\int_0^{t-x}\int_{x}^{t-\alpha}K(t-(\alpha+l))w(x,l)\,dl
f(\alpha)\,d\alpha=[s=t-\alpha]\\
-\int_{-x}^0\int_0^{x+\alpha}K(t-s)w(x,s-\alpha)f(\alpha)dsd\alpha+\\
\int_t^x\int_x^sK(s-l)w(x,l)\,dlf(t-s)\,ds\label{e5}
\end{multline}
Plugging (\ref{e1}), (\ref{e2}), (\ref{e4}), (\ref{e5}) into
(\ref{e3}) and counting the arbitrariness of $f$ and that $f$ vanishes if $\alpha<0$, one get the first
equation in (\ref{Goursat}). Taking $t=x$ in (\ref{e3}) and
choosing $f(0)\not= 0$ implies the second equality in
(\ref{Goursat})

\end{proof}

The space $\mathcal{F}^T:=L_2(0,T)$ with the standard inner
product is said to ba an \emph{outer space} of system
(\ref{dyn_syst1}). The function $u^f(\cdot,t)$ is a state of the
system (\ref{dyn_syst1}) at the moment $t$. The space of states
$\mathcal{H}^T:=L_2(0,T)$ with the standard product is an
\emph{inner space}. The control operator $W^T:\mathcal{F}^T\mapsto
\mathcal{H}^T$ is introduced by the rule:
\begin{equation*}
W^Tf=u^f(x,T).
\end{equation*}
The representation (\ref{u_repr}) implies the \emph{boundary
controllability} of the dynamical system, which is equivalent to
the following
\begin{proposition}
\label{Prop1} The operator $W^T$ is an isomorphism
\end{proposition}
\begin{proof}
In accordance with (\ref{u_repr}) the operator $W^T$ has the
representation
\begin{equation}
\left(W^Tf\right)(x)=f(T-x)+\int_x^Tw(x,s)f(T-s)\,ds.\label{Wtf}
\end{equation}
Thus the boundary controllability is equivalent to the solvability
of the following problem: for prescribed $a\in \mathcal{H}^T$, to
find $f$ such that
\begin{equation}
\label{Volt}
W^Tf=a.
\end{equation}
Then the result of the proposition follows from the fact that
(\ref{Volt}) is a Volterra type equation of the second kind.
\end{proof}

The connecting operator $C^T:\mathcal{F}^T\mapsto \mathcal{F}^T$
is introduced by the rule
\begin{equation}
\label{CT_def}
\left(C^Tf,g\right)_{\mathcal{F}^T}=\left(u^f(\cdot,T),u^g(\cdot,T)\right)_{\mathcal{H}^T}
\end{equation}
Due to the Proposition \ref{Prop1}, $C^T$ is an isomorphism in
$\mathcal{F}^T$. We introduce the Blagoveschenskii function by the
rule
\begin{equation*}
\psi(t,s)=\left(u^f(\cdot,t),u^g(\cdot,s)\right)_{\mathcal{H}^T}
\end{equation*}
The following statement is crucial in the BC method:
\begin{proposition}\label{CtProp}
The operator $C^T$ is determined by inverse data $R^{2T}$ and the function
$K$.
\end{proposition}
\begin{proof}
For $f,g\in C_0^\infty(0,T)$ we evaluate:
\begin{align*}
\psi_{tt}(t,s)=&\int_0^Tu^f_{tt}(x,t)u^g(x,s)\,dx\\
=&\int_0^T\left(u^f_{xx}(x,t)-q(x)u^f(x,t)-\int_0^tK(t-\tau)u^f(x,\tau)\,
d\tau\right)u^g(x,s)\,dx\\
=&\int_0^T u^f(x,t)u_{xx}^g(x,s)-q(x)u^f(x,t)u^g(x,s)\,dx-
u^f_x(0,t)u^g(0,s)\\
&+u^f(0,t)u^g_x(0,s)-\int_0^tK(t-\tau)\int_0^Tu^f(x,\tau)u^g(x,s)\,dx\,d\tau\\
=&\int_0^T
u^f(x,t)u_{xx}^g(x,s)-q(x)u^f(x,t)u^g(x,s)\,dx\\
&-\left(R^Tf\right)(t)g(s)+f(t)\left(R^Tg\right)(s)-\int_0^tK(t-\tau)\psi(\tau,s)\,d\tau.
\end{align*}

Note, that  the terms $u_x^f(T,t)u^g(T,s)$ and $u^f(T,t)u_x^g(T,s)$ are both zero because
of finiteness of $f$ and $g$ and representation formula (\ref{Wtf}).

Evaluating
\begin{align*}
\psi_{ss}(t,s)=&\int_0^Tu^f(x,t)u^g_{ss}(x,s)\,dx\\
=&\int_0^T
u^f(x,t)u_{xx}^g(x,s)-q(x)u^f(x,t)u^g(x,s)\,dx-\int_0^sK(s-\alpha)\psi(t,\alpha)\,d\alpha
\end{align*}
we see that the function $\psi(t,s)$ satisfy the following initial
boundary value problem for the wave equation
\begin{equation}
    \label{eqCT}
    \begin{cases}
        \psi_{tt}(t,s)-\psi_{ss}(t,s)-\int_0^tK(t-\tau)\psi(\tau,s)\,d\tau+\int_0^sK(s-\alpha)\psi(t,\alpha)\,d\alpha\\
        =\left(Rf\right)(t)g(s)-f(t)\left(Rg\right)(s),\\
        \psi(0,s)=\psi(t,0)=\psi_s(0,s)=\psi_t(t,0),
    \end{cases}
\end{equation}
the existence of the solution to (\ref{eqCT}) can be obtained by
the reducing (\ref{eqCT}) to an integral equation, see
\cite{Pan1,Pan2}.

We are left to observe that
\begin{equation*}
\left(C^Tf,g\right)=\psi(T,T).
\end{equation*}

\end{proof}

The operator $C^T$ being bounded operator in $\mathcal{F}^T$ is
given by its kernel $c(t,s)$:
\begin{equation*}
\left(C^Tf\right)(t)=\int_0^Tc(t,s)f(s)\,ds.
\end{equation*}
Since we know the quadratic form (\ref{CT_def}) for all $f,g\in
C_0^\infty(0,T)$, we know the kernel $c(t,s)$.


\section{Gelfand--Levitan equations.}

We introduce the notations:
\begin{eqnarray*}
J_T:\mathcal{F}^T\mapsto \mathcal{F}^T,\quad
\left(J_Tf\right)(t)=f(T-t),\\
W^T=\left(I+M\right)J_T,\quad
\left(Mf\right)(x)=\int_x^Tw(x,s)f(s)\,ds.
\end{eqnarray*}
Then for the inverse operator we have
\begin{eqnarray*}
\left(W^T\right)^{-1}=J_T\left(I+L\right),\quad
\left(Lf\right)(x)=\int_x^Tz(x,s)f(s)\,ds.
\end{eqnarray*}
The equality $(I+M)(I+L)=I$ implies the condition for the kernels
on the diagonal:
\begin{equation*}
z(x,x)=-w(x,x)=\frac{1}{2}\int_0^xq(\alpha)\,d\alpha.
\end{equation*}

For arbitrary $f,g\in \mathcal{F}^T$, by the definition of $C^T$
we have:
\begin{equation}
\label{C_T_R1}
(C^Tf,g)_{\mathcal{F}^T}=(W^Tf,W^Tg)_{\mathcal{H}^T}.
\end{equation}
Let us put $f=(W^T)^{-1}a$, $g=(W^T)^{-1}b$, $a,b\in
\mathcal{H}^T$ and rewrite (\ref{C_T_R1}) as
\begin{equation}
\label{C_T_R2}
(C^TJ_T(I+L)a,J_T(I+L)b)_{\mathcal{F}^T}=(a,b)_{\mathcal{H}^T},
\end{equation}
Since the above equality holds for all $a,b\in \mathcal{H}^T$,
upon introducing the notation $\widetilde C^T:=J_TC^TJ_T$ we see
that (\ref{C_T_R2}) leads to the following operator equation
\begin{equation}
\label{Oper_eqn_inv} (I+L)^*\widetilde C^T(I+L)=I.
\end{equation}
We introduce the notations:
\begin{equation}\label{defC_t}
C_T:=\widetilde C^T-I,\quad \text{and }
\left(C_Tf\right)(t)=\int_0^Tc_T(t,s)f(s)\,ds,
\end{equation}
and rewrite (\ref{Oper_eqn_inv}) as
\begin{equation}
\label{Oper_eqn_1} L^*+(I+L^*)(L+C_T+C_TL)=0.
\end{equation}
Note that the operator $L^*$ has the form:
\begin{equation*}
\left(L^*a\right)(t)=\int_0^tz(x,t)a(x)\,dx
\end{equation*}

The function $z(x,s)$ was defined  for $0\leqslant x\leqslant
s\leqslant T$, let us continue it by zero in the domain
$s<x\leqslant T$ and introduce the function $\phi_x(s)$, $x,s\in
[0,T]$ by the rule
\begin{equation*}
\phi_x(s)=z(x,s)+c_T(x,s)+\int_0^Tc_T(s,\tau)z(x,\tau)\,d\tau.
\end{equation*}
The equality (\ref{Oper_eqn_1}) implies that
\begin{equation*}
z(s,x)+\phi_x(s)+\int_0^T z(s,\tau)\phi_x(\tau)\,d\tau=0,\quad
x,s\in (0,T).
\end{equation*}
Since $z(s,x)=0$ for $0<x<s<T$, we obtain that
\begin{equation*}
\phi_x(s)+\int_s^T z(s,\tau)\phi_x(\tau)\,d\tau=0,\quad 0<x<s<T.
\end{equation*}
Rewriting this equation as
\begin{equation*}
((I+L)\phi_x(\cdot))(t)=0,\quad 0<x<s<T,
\end{equation*}
and taking into account the invertibility of $I+L$ (which follows
from Proposition \ref{Prop1}), we get that
\begin{equation}
\label{G_L}
\phi_x(s)=z(x,s)+c_T(x,s)+\int_x^Tc_T(s,\tau)z(x,\tau)\,d\tau=0,\quad
0<x<s<T.
\end{equation}
Let us formulate this result as
\begin{theorem}
\label{Th_GL} The kernel of operator   $L$
satisfies the following integral equation
\begin{equation}
\label{G_LN}
z(x,s)+c_T(x,s)+\int_x^Tc_T(s,\tau)z(x,\tau)\,d\tau=0,\quad
0<x<s<T.
\end{equation}
where $c_T$ is defined by dynamic inverse data, see  Proposition \ref{CtProp} and (\ref{defC_t}).
\end{theorem}
Solving the equation (\ref{G_LN}) for all $x \in (0,T)$ we can
recover the potential using
\begin{equation*}
q(x)=2\frac{d}{dx}z(x,x).
\end{equation*}

Summing up, note the following:
\begin{itemize}
\item[1)] Using a purely dynamic approach makes it possible to
consider the inverse problem  for the system (\ref {dyn_syst1})
for a more general class of potentials $q\in
L_{1,\,loc}(\mathbb{R}_+)$ see \cite{A2} and kernels $K\in
L_{2,\,loc}(\mathbb{R}_+)$.

\item[2)] The presence of a dynamic representation of the solution
$u^f$ (\ref{u_repr}) (Duhamel representation),  gives the
potential opportunity to use advanced techniques of the Boundary
control method and to study the problem of the characterization of
dynamic inverse data in the spirit of \cite{1D-BC, BM}.

\item[3)] Analysis of the solution of the forward problem (\ref{dyn_syst1}) for general coefficients, analysis of the solution of the equation for the Blagoveshchensky function (\ref{eqCT}), as well as the problem of characterizing inverse data will be the subject of study in the following publications.

\end{itemize}

\noindent{\bf Acknowledgments}

A. S. Mikhaylov and V. S. Mikhaylov were partly supported by RFBR
18-01-00269 , and by Volkswagen Foundation project "From Modeling
and Analysis to Approximation".


\begin{thebibliography}{99}

\bibitem{A1} {S.A. Avdonin, S.A. Ivanov, J.-M. Wang.} \textit{Inverse problems for
the heat equation with memory.} Inverse Problems and Imaging, 13,
no. 1, 31-38, 2019.


\bibitem{A2} {S.A. Avdonin, V.S. Mikhaylov.} \textit{The boundary control approach
to inverse spectral theory.} Inverse Problems, 26, no. 4, 045009,
2010.


\bibitem{A3} {S.A. Avdonin, L. Pandolfi.} \textit{A linear algorithm for the
identification of a weakly singular relaxation kernel using two
boundary measurements.} Journal of Inverse and Ill-Posed Problems
26, no. 2, 299-310, 2018.

\bibitem{B07}
{M.I. Belishev}, \textit{Recent progress in the boundary control
method}, Inverse Problems, 23, no.5, R1-R67, 2007.

\bibitem{1D-BC}
{M.I.Belishev}. \textit{Boundary control and inverse problems:
1-dimensional variant of the BC-method.}
\newblock{\em J. Math. Sciences}, 155, no. 3, 343--379, 2008.

\bibitem{B17}
{M.I.Belishev.} \textit{Boundary control and tomography of
Riemannian manifolds (the BC-method).} Russian Mathematical
Surveys Volume, 72, no. 4, 581-644, 2017.

\bibitem{BM_1}
\textsc{M.I.Belishev and V.S.Mikhailov}. Unified approach to
classical equations of inverse problem theory. \textit{Journal of
    Inverse and Ill-Posed Problems}, 20, no 4, 461--488, 2012.

\bibitem{BM}
{M. I. Belishev, V. S. Mikhaylov.}
\textit{Inverse problem for
    one-dimensional dynamical Dirac system (BC-method)}. Inverse
Problems, 30, no. 12, 2014.


\bibitem{Pan1}  {L. Pandolfi.} \textit{Dynamical identification
of a space varying coefficient in a system with persistent
memory.} Journal of Mathematical Analysis and Applications, 465,
no. 1, 140-158, 2018.


\bibitem{Pan2} {L. Pandolfi}
\textit{Identification of a space varying coefficient of a linear
viscoelastic string of Maxwell-Boltzman type.}
https://arxiv.org/abs/1702.06770


\end{thebibliography}
\end{document}